\documentclass{article}
\usepackage{nips2005e,times}
\usepackage{amssymb}
\usepackage{epsfig}

\newcommand{\mb}[1]{ \mbox{\boldmath$#1$} }
\newcommand{\x}{\mb{x}}
\newcommand{\y}{\mb{y}}
\newcommand{\norm}[1]{\|#1\|}
\newcommand{\ve}{\varepsilon}
\newcommand{\p}{\partial}

\title{Diffusion Maps, Spectral Clustering and Eigenfunctions
of Fokker-Planck Operators}

\author{
Boaz Nadler\thanks{{Corresponding author. http://$\sim$patheon.yale.edu/$\sim$bn29/}}\
\ \  \ \ \ St\'{e}phane Lafon\ \ \ \ \
 Ronald R. Coifman\\
Department of Mathematics, Yale University,
New Haven, CT 06520. \\
\texttt{\{boaz.nadler,stephane.lafon,ronald.coifman\}@yale.edu} 
\And
Ioannis G. Kevrekidis \\
Department of Chemical Engineering and Program in Applied Mathematics\\
Princeton University,
Princeton, NJ 08544 \\
\texttt{yannis@princeton.edu} 
}

\begin{document}

\maketitle

\begin{abstract}
This paper presents a diffusion based probabilistic interpretation of spectral clustering and
dimensionality reduction algorithms that use the eigenvectors of the
normalized graph Laplacian. Given the pairwise adjacency matrix of all points,
we define a diffusion distance between
any two data points and show that the low dimensional
representation of the data by the first few eigenvectors of the corresponding
Markov matrix
is optimal under a certain mean squared error criterion.
Furthermore, assuming that data points are random samples from
a density $p(\x) = e^{-U(\x)}$ we identify these eigenvectors as discrete approximations of
eigenfunctions of a Fokker-Planck operator in a potential $2U(\x)$ with reflecting boundary conditions.
Finally, applying known results regarding the eigenvalues and
eigenfunctions of the continuous Fokker-Planck operator, we provide a mathematical
justification for the success of spectral clustering and dimensional reduction algorithms
based on these first few eigenvectors. This analysis elucidates, in terms of the characteristics
of diffusion processes, many empirical findings regarding spectral clustering algorithms.
\end{abstract}


\section{Introduction}

Clustering and low dimensional representation of high dimensional
data are important problems in many diverse fields. In recent years
various spectral methods to perform these tasks, based on the
eigenvectors of adjacency matrices of graphs on the data have been
developed, see for example \cite{Scholkopf}-\cite{PNAS} and references therein. In the simplest version, known as
the normalized graph Laplacian, given $n$ data points
$\{\x_i\}_{i=1}^n$ where each $\x_i\in\mathbb{R}^p$, we define a
pairwise similarity matrix between points, for example using a Gaussian kernel with width $\ve$,
\begin{equation}
L_{i,j} = k(\x_i,\x_j) = \exp\left(-\frac{\norm{\x_i - \x_j}^2}{2\ve}\right)
    \label{L_ij}
\end{equation}
and a diagonal normalization matrix $D_{i,i} = \sum_j L_{i,j}$.
Many works propose to use the first few eigenvectors of
the normalized eigenvalue problem $L\phi=\lambda D \phi$, or equivalently of
the matrix $M = D^{-1}L$, either as a low dimensional representation of data
or as good coordinates for clustering purposes. Although eq. (1) is based on a Gaussian kernel,
other kernels are possible.
While for actual datasets the choice of a kernel $k(\x_i,\x_j)$ is crucial,  it does not qualitatively change
our asymptotic analysis \cite{Acha}.

The use of the first few eigenvectors of $M$ as good coordinates is typically justified
with heuristic arguments or as a relaxation of a discrete clustering problem \cite{Shi_Malik}.
In \cite{Belkin02,Belkin03} Belkin and Niyogi showed that when data is uniformly
sampled from a low dimensional manifold of $\mathbb{R}^p$ the first few eigenvectors of $M$
are discrete approximations of the eigenfunctions of the Laplace-Beltrami operator on the manifold,
thus providing a mathematical justification for their use in this case.
A different theoretical analysis of the eigenvectors of the matrix $M$, based on the
fact that $M$ is a stochastic matrix and thus represents a random
walk on the graph has been described by Meil\v{a} and Shi
\cite{Meila}, who considered the case of piecewise constant eigenvectors for
specific lumpable matrix structures. Another notable work that considered the random walk
aspects of spectral clustering is \cite{Saerens,Fouss}, where the authors suggest
clustering based on the average commute time between points.

In this paper we provide a unified probabilistic framework which
combines these results and extends them in two different directions.
First, in section \ref{s:Diffusion_distance} we define a distance
function between any two points based on the random walk on the
graph, which we naturally denote the {\em diffusion distance}. We
then show that the low dimensional description of the data by
the first few eigenvectors, denoted as the {\em diffusion map}, is optimal under a mean squared error
criterion based on this distance. In section \ref{s:FP} we
consider a statistical model, in which data points are iid random samples from a probability
density $p(\x)$ in a smooth bounded domain $\Omega\subset\mathbb{R}^p$ and analyze the asymptotics of the eigenvectors
as the number of data points tends to infinity. This analysis shows
that the eigenvectors of the finite matrix $M$ are discrete approximations of the
eigenfunctions of a Fokker-Planck (FP) operator with reflecting
boundary conditions. This observation, coupled with known results
regarding the eigenvalues and eigenfunctions of
the FP operator provide new insights into the properties of these
eigenvectors and on the performance of spectral clustering algorithms, as described
in section \ref{s:Spectral_Gap}.

\section{Diffusion Distances and Diffusion Maps}\label{s:Diffusion_distance}

The starting point of our analysis, as also noted in other works, is the observation that the
matrix $M$ is adjoint to a symmetric matrix
\begin{equation}
M_{s} = D^{1/2} M D^{-1/2}
\end{equation}
Thus, $M$ and $M_s$ share the same eigenvalues. Moreover, since
$M_s$ is symmetric it is diagonalizable and has a set of $n$
eigenvalues $\{\lambda_j\}_{j=0}^{n-1}$
whose corresponding eigenvectors $\{\mb{v}_{j}\}$ form an
orthonormal basis of $\mathbb{R}^n$. The left and right eigenvectors
of $M$, denoted $\phi_j$ and $\psi_j$ are related to those of $M_s$
according to
\begin{equation}
\phi_j = \mb{v}_j D^{1/2}, \qquad
\psi_j = \mb{v}_j D^{-1/2}\label{phi_psi}
\end{equation}
Since the eigenvectors $\mb{v}_j$ are orthonormal under the standard dot
product in $\mathbb{R}^n$, it follows that the vectors $\phi_j$ and $\psi_k$
are bi-orthonormal
\begin{equation}
  \langle \phi_i,\psi_j\rangle = \delta_{i,j}
    \label{bi-orthonormal}
  \end{equation}
where $\langle \mb{u},\mb{v}\rangle$ is the standard dot product between two vectors
in $\mathbb{R}^n$.
We now utilize the fact that by construction $M$ is a stochastic matrix with
all row sums equal to one, and can thus be interpreted as defining a random walk
on the graph. Under this view, $M_{i,j}$
denotes the transition probability from the point $\x_i$ to the point $\x_j$
in one time step. Furthermore, based on the similarity of the Gaussian kernel
(\ref{L_ij}) to the fundamental solution of the heat equation, we define our time step
as $\Delta t = \ve$. Therefore,
\begin{equation}
  \Pr\{\x(t+\ve) = \x_j\,|\,\x(t) = \x_i\} = M_{i,j}
\end{equation}
Note that $\ve$ has therefore a dual interpretation in this framework. The first is
that $\ve$ is the (squared) radius of the neighborhood used to infer local geometric and density information
for the construction of the adjacency matrix, while the second is that $\ve$ is
the discrete time step at which the random walk jumps from point to point.

We denote by $p(t,\y|\x)$ the probability distribution of a random walk landing at location
$\y$ at time $t$, given a starting location $\x$ at time $t=0$. For $t=k\,\ve$,
$p(t,\y|\x_i) = \mb{e}_i M^k$, where $\mb{e}_i$ is a row vector of zeros with a single one
at the $i$-th coordinate.
For $\ve$ large enough, all points in the
graph are connected so that $M$ has a unique eigenvalue equal to 1.
The other eigenvalues form a non-increasing sequence of non-negative numbers:
$\lambda_0=1 > \lambda_1\geq \lambda_2\geq\ldots\geq
\lambda_{n-1}\geq 0$. Then, regardless of the initial starting point $\x$,
\begin{equation}
  \lim_{t\to\infty} p(t,\y|\x) = \phi_0(\y)
\end{equation}
where $\phi_0$ is the left eigenvalue of $M$ with eigenvalue $\lambda_0 = 1$, explicitly
given by
\begin{equation}
\phi_0(\x_i) = \frac{{D_{i,i} }}
{{\sum\limits_j {D_{j,j} } }}
\end{equation}
This eigenvector also has a dual interpretation. The first is that $\phi_0$ is
the stationary probability distribution on the graph, while the second is that $\phi_0(\x)$ is the density
estimate at the point $\x$. Note that for a general shift invariant kernel $K(\x-\y)$
and for the Gaussian kernel in particular, $\phi_0$ is simply the well known Parzen window
density estimator.

For any finite time $t$, we decompose the probability distribution in the
eigenbasis $\{\phi_j\}$
\begin{equation}
  p(t,\y|\x) = \phi_0(\y) + \sum_{j\geq 1} a_j(\x) \lambda_j^t \phi_j(\y)
    \label{p_t_y}
\end{equation}
where the coefficients $a_j$ depend on the initial location $\x$.
Using the bi-orthonormality condition (\ref{bi-orthonormal}) gives $a_j(\x)=\psi_j(\x)$,
with $a_0(\x)=\psi_0(\x) = 1$ already implicit in (\ref{p_t_y}).


Given the definition of the random walk on the graph it is only natural to quantify the
similarity between any two points based on the evolution of their
probability distributions. Specifically, we consider the following
distance measure at time $t$,
\begin{eqnarray}
  D^2_t(\x_0,\x_1) &=& \norm{p(t,\y|\x_0) - p(t,\y|\x_1)}_w^2 \label{diffusion_distance}\\
  &=& \sum_{\y} (p(t,\y|\x_0) - p(t,\y|\x_1))^2 w(\y) \nonumber
\end{eqnarray}
with the specific choice $w(\y) = 1/\phi_0(\y)$ for the weight function, which
takes into account the (empirical) local density of the points.

Since this distance depends on the random walk on the graph, we quite naturally
denote it as the {\em diffusion distance} at time $t$. We also denote the mapping
between the original space and the first $k$ eigenvectors as the {\em diffusion map}
\begin{equation}
  \Psi_t(\x) = \left(\lambda_1^t \psi_1(\x),\lambda_2^t \psi_2(\x),\ldots,\lambda_k^t \psi_k(\x)\right)
\end{equation}
The following theorem relates the diffusion distance and the diffusion map.

{\bf Theorem:} The diffusion distance (\ref{diffusion_distance}) is equal to Euclidean distance in the diffusion map space
with all $(n-1)$ eigenvectors.
\begin{equation}
D_t^2(\x_0,\x_1) = \sum_{j\geq 1} \lambda_j^{2t}
\left(\psi_j(\x_0) - \psi_j(\x_1)\right)^2 = \norm{\Psi_t(\x_0)-\Psi_t(\x_1)}^2
\end{equation}
{\bf Proof:} Combining (\ref{p_t_y}) and (\ref{diffusion_distance}) gives
\begin{equation}
  D_t^2(\x_0,\x_1) = \sum_{\y} \Big(
  \sum_j \lambda_j^t(\psi_j(\x_0)-\psi_j(\x_1))\phi_j(\y)\Big)^2 1/\phi_0(\y)
\end{equation}
Expanding the brackets, exchanging the order of summation and using relations (\ref{phi_psi})
and (\ref{bi-orthonormal})
between $\phi_j$ and $\psi_j$ yields the required result. Note that the weight factor $1/\phi_0$ is
essential for the theorem to hold. \hfill$\Box$.

This theorem provides a justification for using
Euclidean distance in the diffusion map space for spectral clustering purposes.
Therefore, geometry in diffusion space is meaningful and can be interpreted in terms of the
Markov chain. In particular, as shown in \cite{Lafon_Lee}, quantizing this diffusion space
is equivalent to lumping the random walk. Moreover,
since in many practical applications the spectrum of the matrix $M$ has
a {\em spectral gap} with only a few eigenvalues close to one and all additional eigenvalues
much smaller than one, the diffusion distance at a large enough time $t$ can be well approximated by
only the first few $k$ eigenvectors $\psi_1(\x),\ldots,\psi_k(\x)$, with a
negligible error of the order of $O((\lambda_{k+1}/\lambda_k)^t)$.
This observation provides a theoretical justification for dimensional reduction with
these eigenvectors.
In addition, the following theorem shows that this $k$-dimensional approximation is {\em optimal}
under a certain mean squared error criterion.

{\bf Theorem:} Out of all $k$-dimensional approximations of the form
\[
\hat{p}(t,\y|\x) = \phi_0(\y) + \sum_{j=1}^k a_j(\x,t) \mb{w}_j(\y)
\]
for the probability distribution at time $t$, the one that minimizes the
mean squared error
\[
\mathbb{E}_{\x}\{\norm{p(t,\y|\x) - \hat{p}(t,\y|\x)}_{w}^2\}
\]
where averaging over initial points $\x$ is with respect to the stationary density $\phi_0(\x)$, is given by
$\mb{w}_j(\y) = \phi_j(\y)$ and $a_j(t,\x) = \lambda_j^t\psi_j(\x)$.
Therefore, the optimal $k$-dimensional
approximation is given by the truncated sum
\begin{equation}
\hat{p}(\y,t|\x) = \phi_0(\y) + \sum_{j=1}^k \lambda_j^t\psi_j(\x)\phi_j(\y)
\end{equation}
{\bf Proof:} The proof is a consequence of a weighted principal component analysis
applied to the matrix $M$, taking into account the biorthogonality of the left and right eigenvectors.

We note that the first few eigenvectors are also optimal under other
criteria, for example for data sampled from a manifold as in \cite{Belkin02},
or for multiclass spectral clustering \cite{Yu}.

\section{The Asymptotics of the Diffusion Map}
\label{s:FP}

The analysis of the previous section provides a mathematical
explanation for the success of the diffusion maps
for dimensionality reduction and spectral clustering. However, it does not provide
any information regarding the structure of the computed
eigenvectors.

To this end, we introduce a statistical model and assume that the data points $\{\x_i\}$
are i.i.d. random samples from a probability density $p(\x)$ confined to a compact
connected subset $\Omega\subset\mathbb{R}^p$ with smooth boundary $\p\Omega$.
Following the statistical physics notation, we write the density in Boltzmann form, $p(\x) = e^{-U(\x)}$,
where $U(\x)$ is the (dimensionless) potential or energy of the configuration $\x$.

As shown in \cite{Acha}, in the limit $n\to\infty$ the random walk on the discrete graph
converges to a random walk on the continuous space $\Omega$. Then, it is possible to define
forward and backward operators $T_f$ and $T_b$ as follows,
\begin{equation}
\begin{array}{*{20}c}
   T_f [\phi ](\x) \hfill &  =  \hfill & \int_{\Omega} M(\x|\y) \phi(\y) p(\y) d\y  \hfill  \\
    & & \\
   T_b [\psi](\x)  \hfill &  =  \hfill & \int_{\Omega}M(\y|\x) \psi(\y) p(\y) d\y  \hfill  \\
 \end{array}
\end{equation}
where $M(\x|\y) = \exp(-\norm{\x-\y}^2/2\ve)/D(\y)$ is the transition probability from $\y$
to $\x$ in time $\ve$, and $D(\y) = \int \exp(-\norm{\x-\y}^2/2\ve)p(\x) d\x$.

The two operators $T_f$ and $T_b$ have probabilistic interpretations. If $\phi(\x)$ is a probability
distribution on the graph at time $t=0$, then $T_f[\phi]$ is the probability distribution
at time $t=\ve$. Similarly, $T_b[\psi](\x)$ is the mean of the function $\psi$ at time $t=\ve$,
for a random walk that started at location $\x$ at time $t=0$.
The operators $T_f$ and $T_b$ are thus the continuous analogues of the left and right multiplication
by the finite matrix $M$.
In \cite{Luxburg} a rigorous mathematical proof of the convergence of the eigenvalues and
eigenvectors of the discrete matrix $M$ to the eigenvalues and eigenfunctions of
the integral operators $T_f$ and $T_b$ is given.

We now take this analysis one step further and consider the limit $\ve\to 0$.
This is possible, since when $n=\infty$ each data point contains an infinite number of nearby neighbors.
In this limit, since $\ve$ also has the interpretation of a time step, the random walk converges to a
diffusion process, whose probability density evolves continuously in time, according to
\begin{equation}
  \frac{\partial p(x,t)}{\partial t}=\lim_{\ve\to 0}\frac{p(x,t+\ve)-p(x,t)}{\ve}=
  \lim_{\ve\to 0} \frac{T_f-I}\ve p(x,t)
\end{equation}
 in which case it is customary to study the infinitesimal generators
(propagators)
\begin{equation}
  {\cal H}_f = \lim_{\ve \to 0}\frac{T_f - I }{\ve},\qquad
    {\cal H}_b = \lim_{\ve \to 0}\frac{T_b - I }{\ve} \label{propagators}
\end{equation}
Clearly, the eigenfunctions of $T_f$ and $T_b$ converge to those of ${\cal H}_f$ and ${\cal H}_b$, respectively.

As shown in \cite{Acha}, the backward generator is given by
the following Fokker-Planck operator
\begin{equation}
  {\cal H}_b \psi = \Delta \psi - 2 \nabla \psi \cdot \nabla U
    \label{Hb}
\end{equation}
which corresponds to a diffusion process in a potential field of $2U(\x)$
\begin{equation}
  \dot{\x}(t) = - \nabla (2U) + \sqrt{2D} \dot{\mb{w}}(t) \label{Langevin}
\end{equation}
where $\mb{w}(t)$ is standard Brownian motion in $p$ dimensions and $D$ is the
diffusion coefficient, equal to one in equation (\ref{Hb}).
The Langevin equation (\ref{Langevin}) is the most common equation to describe
stochastic dynamical systems in physics, chemistry and biology \cite{Gardiner,Risken}.
As such, its characteristics as well as those of the corresponding FP equation
have been extensively studied, see
\cite{Gardiner}-\cite{Eckhoff} and many others.
The term $\nabla \psi \cdot \nabla U$ in (\ref{Hb}) is interpreted
as a {\em drift} term towards low energy (high-density) regions, and as discussed
in the next section, may play a crucial part in the definition of clusters.

Note that when data is uniformly sampled from $\Omega$, $\nabla U=0$ so the drift term
vanishes and we recover the Laplace-Beltrami operator on $\Omega$.
The convergence (in probability) of the discrete matrix $M$
to this operator has been recently rigorously proven in \cite{Belkin05}. However,
the important issue of boundary conditions was not considered.

Since (\ref{Hb}) is defined in the bounded domain $\Omega$,
the eigenvalues and eigenfunctions of ${\cal H}_b$ depend on the boundary
conditions imposed on $\p\Omega$. As shown in \cite{Lafon}, in the limit $\ve\to 0$, the random walk satisfies reflecting
boundary conditions on $\p\Omega$, which translate into
\begin{equation}
  \frac{\p\psi(\x)}{\p \mb{n}}\Big|_{\p\Omega} = 0
\end{equation}
where $\mb{n}$ is a unit normal vector at the point $\x\in\p\Omega$.

To conclude, the left and right eigenvectors of the finite matrix $M$ can be viewed as discrete approximations
to those of the operators $T_f$ and $T_b$, which in turn can be viewed as approximations to those of ${\cal H}_f$
and ${\cal H}_b$. Therefore, if there are enough data points for
accurate statistical sampling, the structure and characteristics of the eigenvalues and
eigenfunctions of ${\cal H}_b$
are similar to the corresponding eigenvalues and discrete eigenvectors of $M$. For convenience, the three different stochastic processes are shown in table 1.

\begin{table}[t]
\caption{Random Walks and Diffusion Processes}
\begin{center}{
\begin{tabular}{|c|c|c|}
\hline
{\bf Case} & {\bf Operator} & {\bf Stochastic Process}\\
\hline
$\ve > 0$  & finite $n\times n$ & R.W. discrete in space  \\
$n<\infty$ & matrix  $M$& discrete in time \\
\hline 
$\ve > 0$  & operators & R.W. continuous in space  \\
$n\to\infty$ & $T_f, T_b$ & discrete in time \\
\hline 
$\ve \to 0$  & infinitesimal  & diffusion process  \\
$n=\infty$ & generator ${\cal H}_f$ & continuous in time \& space\\
\hline
\end{tabular}
}
\end{center}
\end{table}

\section{Fokker-Planck eigenfunctions and spectral clustering}\label{s:Spectral_Gap}

According to (\ref{propagators}), if $\lambda_\ve$ is an eigenvalue of the matrix $M$ or of the integral
operator $T_b$ based on a kernel with parameter $\ve$, then the corresponding eigenvalue
of ${\cal H}_b$ is $\mu \approx (\lambda_\ve - 1)/\ve$. Therefore the largest eigenvalues of $M$
correspond to the smallest eigenvalues of ${\cal H}_b$.
These eigenvalues and their corresponding eigenfunctions have been extensively studied in
the literature under various settings. In general, the eigenvalues and eigenfunctions depend
both on the geometry of the domain $\Omega$ and on the profile of the potential $U(\x)$.
For clarity and due to lack of space we briefly analyze here two extreme cases. In the first case
$\Omega=\mathbb{R}^p$ so geometry plays no role, while in the second $U(\x)=const$
so density plays no role. Yet we show that in both cases there can still be well defined clusters,
with the unifying probabilistic concept being that the mean exit time from one cluster to another
is much larger than the characteristic equilibration time inside each cluster.

{\bf Case I:} Consider diffusion in a smooth potential field $U(\x)$ in $\Omega=\mathbb{R}^p$, where
$U$ has a few local minima, and $U(\x)\to\infty$ as $\norm{\x}\to\infty$ fast enough
so that $\int e^{-U} d\x = 1 <\infty$. Each such local minimum thus defines
a metastable state, with transitions between metastable states being relatively rare
events, depending on the barrier heights separating them.
As shown in \cite{Matkowsky,{Eckhoff}} (and in many other works) there is an intimate connection between
the smallest eigenvalues
of ${\cal H}_b$ and mean exit times out of these metastable states. Specifically, in the
asymptotic limit of small noise $D\ll 1$, exit times are exponentially distributed and the first
non-trivial eigenvalue (after $\mu_0=0$) is given by $\mu_1 = 1/\bar{\tau}$ where
$\bar{\tau}$ is the mean exit time to overcome the highest potential barrier on the way to the deepest
potential well. For the case of two potential wells, for example,
the corresponding eigenfunction is roughly constant in each well with a sharp transition
near the saddle point between the wells.
In general, in the case of $k$ local minima there are asymptotically only $k$ eigenvalues very close to zero.
Apart from $\mu_0=0$, each of the other $k-1$ eigenvalues corresponds to the mean exit time from
one of the wells into the deepest one, with the corresponding eigenfunctions being almost constant in each well.
Therefore, for a finite dataset the presence of only $k$ eigenvalues close to 1 with a
{\em spectral gap}, e.g. a large difference between $\lambda_k$ and $\lambda_{k+1}$ is indicative
of $k$ well defined {\em global} clusters. In figure 1 (left) an example of this case is
shown, where $p(\x)$ is the sum of two well separated Gaussian clouds leading to a double well potential.
Indeed there are only two eigenvalues close or equal to 1 with a distinct spectral gap and the first
eigenfunction being almost piecewise constant in each well.

In stochastic dynamical systems a spectral gap corresponds to a separation of time scales between
long transition times from one well or metastable state to another as compared to short
equilibration times inside each well. Therefore, clustering and identification of metastable
states are very similar tasks, and not surprisingly algorithms similar to
the normalized graph Laplacian have been independently
developed in the literature \cite{Schutte}.

The above mentioned results are asymptotic in the small noise limit. In practical datasets, there can be
clusters of different scales, where a global analysis with a single $\ve$ is
not suitable. As an example consider the second dataset in figure 1, with three clusters.
While the first eigenvector distinguishes between the large cluster and the two smaller ones,
the second eigenvector captures
the equilibration inside the large cluster instead of further distinguishing the two small clusters.
While a theoretical explanation is beyond the scope of this paper, a possible solution is to choose
a location dependent $\ve$, as proposed in \cite{Zelnik}.

{\bf Case II:} Consider a uniform density in a region
$\Omega\subset\mathbb{R}^3$ composed of two large containers connected by a narrow circular tube,
as in the top right frame in figure 1. In this case $U(\x) = const$,
so the second term in (\ref{Hb}) vanishes. As shown in \cite{Singer}, the second eigenvalue
of the FP operator is extremely small, of the order of $a/V$ where $a$ is the radius of
the connecting tube and $V$ is the volume of the containers, thus showing an interesting connection
to the Cheeger constant on graphs.
The corresponding eigenfunction is
almost piecewise constant in each container with a sharp transition in the connecting tube. Even though in this case the density
is uniform, there still is a spectral gap with two well defined clusters (the two containers),
defined entirely by the geometry of $\Omega$. An example of such a case and the results of the
diffusion map are shown in figure 1 (right).

In summary the eigenfunctions and eigenvalues of the FP operator, and thus
of the corresponding finite Markov matrix,
depend on both geometry and density. The diffusion distance and its close relation
to mean exit times between different clusters is the quantity that incorporates these
two features. This provides novel insight
into spectral clustering algorithms, as well as a theoretical justification
for the algorithm in \cite{Fouss}, which defines clusters according
to mean travel times between points on the graph. Finally, for dynamical systems these eigenvectors
can be used to design better search and data collection protocols \cite{Yannis}.


\begin{figure}[t]
\mbox{
\begin{minipage}[t] {\textwidth}
\begin{center}
\begin{tabular}{c}
\psfig{figure=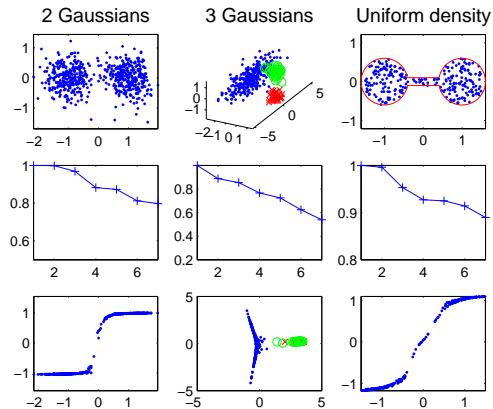,width=6.5cm}
\end{tabular}
\end{center}
\end{minipage}
}
\caption{Diffusion map results on different datasets. Top - the datasets. Middle - the eigenvalues. Bottom -
the first eigenvector vs. $\x_1$ or the first and second eigenvectors for the case of three Gaussians.}
\label{f:cls1}
\end{figure}

\subsubsection*{Acknowledgments}
The authors thank Mikhail Belkin and Partha Niyogi for interesting
discussions. This work was partially supported by DARPA.


\small{

\end{document}
\begin{thebibliography}{*99}

\bibitem{Scholkopf}
B. Sch\"{o}lkopf, A. Smola and K.R. M\"{u}ller. Nonlinear component analysis
as a kernel eigenvalue problem, {\em Neural Computation} 10, 1998.

\bibitem{Weiss} Y. Weiss. Segmentation using eigenvectors: a
unifying view. {\em ICCV} 1999.

\bibitem{Shi_Malik}
J. Shi and J. Malik. Normalized cuts and image segmentation,
{\em PAMI}, Vol. 22, 2000.

\bibitem{Belkin02}
M. Belkin and P. Niyogi. Laplacian eigenmaps and spectral techniques
for embedding and clustering, NIPS Vol. 14, 2002.

\bibitem{Belkin03}
M. Belkin and P. Niyogi. Laplacian eigenmaps for dimensionality reduction
and data representation, {\em Neural Computation} 15:1373-1396, 2003.

\bibitem{Ng}
A.Y. Ng, M. Jordan and Y. Weiss. On spectral clustering, analysis and
an algorithm, NIPS Vol. 14, 2002.

\bibitem{Saerens}
M. Saerens, F. Fouss, L. Yen and P. Dupont,
The principal component analysis of a graph and its relationships to
spectral clustering.
Proceedings of the 15th European Conference on Machine Learning (ECML 2004),
Lecture Notes in Artificial Intelligence, Vol. 3201, 
2004.

\bibitem{Lafon}
R.R. Coifman, S. Lafon,
Diffusion Maps, to appear in Appl. Comp. Harm. Anal.

\bibitem{PNAS} R.R. Coifman \& {\em al.}, 
Geometric diffusion as a tool for harmonic
analysis and structure definition of data, parts I and II, {\em Proc.
Nat. Acad. Sci.}, in press.

\bibitem{Acha} B. Nadler, S. Lafon, R.R. Coifman, I. G. Kevrekidis, Diffusion
maps, spectral clustering, and the reaction coordinates of
dynamical systems, to appear in Appl. Comp. Harm. Anal., available at
http://arxiv.org/abs/math.NA/0503445.


\bibitem{Meila} M. Meila, J. Shi.
A random walks view of spectral segmentation,
{\em AI and Statistics}, 2001.

\bibitem{Fouss}
L. Yen L., Vanvyve D., Wouters F., Fouss F., Verleysen M. and Saerens M. ,
Clustering using a random-walk based distance measure.
Proceedings of the 13th Symposium on Artificial Neural Networks (ESANN 2005), pp 317-324 (2005).

\bibitem{Yu}
S. Yu and J. Shi. Multiclass spectral clustering. ICCV 2003.

\bibitem{Luxburg}
U. von Luxburg, O. Bousquet, M. Belkin,
On the convergence of spectral clustering on random samples: the
normalized case, NIPS, 2004.

\bibitem{Lafon_Lee}
S. Lafon, A.B. Lee,
Diffusion maps: A unified framework for dimension reduction,
data partitioning and graph subsampling, submitted.

\bibitem{Gardiner}
C.W. Gardiner, {\em Handbook of stochastic methods}, third edition, Springer NY, 2004.

\bibitem{Risken}
H. Risken, {\em The Fokker Planck equation}, 2nd edition, Springer NY, 1999.

\bibitem{Matkowsky}
B.J. Matkowsky and Z. Schuss, Eigenvalues of the Fokker-Planck
operator and the approach to equilibrium for diffusions in
potential fields, {\em SIAM J. App. Math.} 40(2):242-254 (1981).


\bibitem{Eckhoff}
M. Eckhoff, Precise asymptotics of small eigenvalues of reversible diffusions in the metastable
regime, {\em Annals of Prob.} 33:244-299, 2005.

\bibitem{Belkin05}
M. Belkin and P. Niyogi, Towards a theoeretical foundation for
Laplacian-based manifold methods, COLT 2005 (to appear).

\bibitem{Schutte}
W. Huisinga, C. Best, R. Roitzsch, C. Sch\"{u}tte, F. Cordes, From simulation data
to conformational ensembles, structure and dynamics based methods, {\em J. Comp. Chem.}
20:1760-74, 1999.


\bibitem{Zelnik}
L. Zelnik-Manor, P. Perona,
Self-Tuning spectral clustering, NIPS, 2004.


\bibitem{Singer}
A. Singer, Z. Schuss, D. Holcman and R.S. Eisenberg, narrow escape, part I, submitted.

\bibitem{Yannis}
I.G. Kevrekidis, C.W. Gear, G. Hummer,  Equation-free:
The computer-aided analysis of comptex multiscale systems,
{\em Aiche J. } 50:1346-1355, 2004.

\end{thebibliography}
